\documentclass[a4paper,11pt]{article}
\usepackage[latin1]{inputenc}
\usepackage{epsfig}
\usepackage{amsmath,amssymb}
\usepackage{amsfonts}
\usepackage{indentfirst}
\date{}

\title{Enumerating permutations avoiding more than three Babson - Steingr\'\i msson patterns}
\author{Antonio Bernini, Elisa Pergola\\
\\
\small{Dipartimento di Sistemi e
Informatica, Università di Firenze}\\
{\small viale G. B. Morgagni 65, 50134 Firenze, Italy}\\
{\small\tt{bernini@dsi.unifi.it; elisa@dsi.unifi.it}}}

\newtheorem{prop}{Proposition}[section]

\newcommand{\qed}{\begin{flushright}$\square$\end{flushright}
}

\begin{document}

\maketitle

\begin{abstract}
Not long ago, Claesson and Mansour proposed some conjectures about
the enumeration of the permutations avoiding more than three
Babson - Steingr\'\i msson patterns (generalized patterns of type
$(1,2)$ or $(2,1)$). The avoidance of one, two or three patterns
has already been considered. Here, the cases of four and five
forbidden patterns are solved and the exact enumeration of the
permutations avoiding them is given, confirming the conjectures of
Claesson and Mansour. The approach we use can be easily extended
to the cases of more than five forbidden patterns.
\end{abstract}

\section{Introduction}

The results of the present paper concern the exact enumeration of
the permutations, according to their length, avoiding any set of
four or five generalized patterns \cite{BS} of type $(1,2)$ or
$(2,1)$. The cases of the permutations avoiding one, two or three
generalized patterns (of the same types) were solved in \cite{C},
\cite{CM} and \cite{BFP}, respectively. In particular, in
\cite{CM} the authors conjectured the plausible sequences
enumerating the permutations of $S_n(P)$, for any set $P$ of three
or more patterns.

In \cite{BFP}, the proofs were substantially conducted by finding
the ECO construction \cite{BDPP} for the permutations avoiding
three generalized patterns of type $(1,2)$ or $(2,1)$, encoding it
with a succession rule and, finally, checking that this one leads
to the enumerating sequence conjectured in \cite{CM}. This
approach could be surely used also for the investigation of the
avoidance of four or five generalized patterns of type $(1,2)$ or
$(2,1)$ and, maybe, it would allow to find same nice and
interesting results: we think that, for instance, in same case new
succession rules for known sequences would appear. Nevertheless,
this approach has just one obstacle: the large number of cases to
consider in order to exhaust all the conjectures in \cite{CM}. The
line we are going to follow (see below) is simple and allows us to
reduce the number of cases to be considered. Most of the results
are summarized in several tables which are presented in Section
\ref{tables}. Really, the paper could appear an easy exercise, but
we believe that it is a valuable contribute to the classification
of permutations avoiding generalized patterns, started with
Claesson, Mansour, Elizalde and Noy \cite{EN}, Kitaev \cite{K}.
Moreover, it can be seen as the continuation of the work started
in \cite{BFP} for the fulfillment of the proofs of the conjectures
presented in \cite{CM}.

\subsection{Preliminaries}

A (classical) \emph{pattern} is a permutation $\sigma\in S_k$ and
a permutation $\pi\in S_n$ \emph{avoids} $\sigma$ if there is no
any subsequence $\pi_{i_1}\pi_{i_2}\ldots\pi_{i_k}$ with $1\leq
i_1<i_2<\ldots<i_k\leq n$ which is order-isomorphic to $\sigma$.
In other word, $\pi$ must contain no subsequences having the
entries in the same relative order of the entries of $\sigma$.
\emph{Generalized patterns} were introduced by Babson and
Steingr\'imsson for the study of the mahonian statistics on
permutations \cite{BS}. They are constructed by inserting one or
more dashes among the elements of a classical pattern (two or more
consecutive dashes are not allowed). For instance, $216-4-53$ is a
generalized pattern of length $6$. The \emph{type}
$(t_1,t_2,\ldots,t_{h+1})$ of a generalized pattern containing $h$
dashes records the number of elements between two dashes (we
suppose a dash at the beginning and at the end of the generalized
pattern, but we omit it): the type of $216-4-53$ is $(3,1,2)$. A
permutation $\pi$ \emph{contains} a generalized pattern $\tau$ if
$\pi$ contains $\tau$ in the classical sense and if any pair of
elements of $\pi$ corresponding to two adjacent elements of $\tau$
(not separated by a dash) are adjacent in $\pi$, too. For
instance, $\pi=153426$ contains $32-14$ in the entries
$\pi_2\pi_3\pi_5\pi_6=5326$ or the pattern $3-214$ in the entries
$\pi_2\pi_4\pi_5\pi_6=5426$. A permutation $\pi$ \emph{avoids} a
generalized pattern $\tau$ if it does not contain $\tau$. If $P$
is a set of generalized patterns, we denote $S_n(P)$ the
permutations of length $n$ of $S$ (symmetric group) avoiding the
patterns of $P$.

In the paper, we are interested to the generalized patterns of
length three, which are of type $(1,2)$ or $(2,1)$ and are those
ones specified in the set
$$
\begin{array}{ll}
\mathcal M =& \{1-23,12-3,1-32,13-2,3-12,31-2,2-13,21-3,\\
&\mbox{  } \mbox{  } 2-31,23-1,3-21,32-1\}.\\
\end{array}
$$
In the sequel, sometimes we can refer to a \emph{generalized
pattern of length three} more concisely with \emph{pattern}.

\bigskip\noindent
If $\pi\in S$, we define its \emph{reverse} and its
\emph{complement} to be the permutations $\pi^r$ and $\pi^c$,
respectively, such that $\pi^r_i=\pi_{n+1-i}$ and
$\pi_i^c=n+1-\pi_i$. We generalize this definition to a
generalized pattern $\tau$ obtaining its reverse $\tau^r$ by
reading $\tau$ from right to left (regarding the dashes as
particular entries) and its complement $\tau^c$ by considering the
complement of $\tau$ regardless of the dashes which are left
unchanged (e.g. if $\tau=216-4-53$, then $\tau^r=35-4-612$ and
$\tau^c=561-3-24$). It is easy to check $\tau^{rc}=\tau^{cr}$. If
$P\subseteq\mathcal M$, the set $\{P,P^r,P^c,P^{rc}\}$ is called
the \emph{symmetry class} of $P$ ($P^r,P^c\ \mbox{and}\ P^{rc}$
contain the reverses, the complements and the reverse-complements
of the patterns specified in $P$, respectively). We have that
$|S_n(P)|$=$|S_n(P^r)|$=$|S_n(P^c)|$=$|S_n(P^{rc})|$ \cite{SS},
therefore we can choose one of the four possible sets as the
\emph{representative} of a symmetry class, as far as the
enumeration of $S(R),\ R\in\{P,P^r,P^c,P^{rc}\}$, is concerned.

\subsection{The strategy}
Looking at the table of \cite{CM} where the authors present their
conjectures, it is possible to note that most of the sequences
enumerating the permutations avoiding four patterns are the same
of those ones enumerating the permutations avoiding three
patterns. A similar fact happens when the forbidden patterns are
four and five. This suggests to use the results for the case of
three forbidden patterns (at our disposal) to deduce the proof of
the conjectures for the case of four forbidden patterns and,
similarly, use the results for the case of four forbidden patterns
to solve the case of five forbidden patterns. Indeed, it is
obvious that $S(p_1,p_2,p_3,p_4)\subseteq
S(p_{i_1},p_{i_2},p_{i_3})$ (with $i_j\in \{1,2,3,4\}$ and $p_l\in
\mathcal M$). If the inverse inclusion can be proved for some
patterns, then the classes $S(p_1,p_2,p_3,p_4)$ and
$S(p_{i_1},p_{i_2},p_{i_3})$ coincide and they are enumerated by
the same sequence (a similar argument can be used for the case of
four and five forbidden patterns).

The following eight propositions are useful to this aim, as well:
each of them proves that if a permutation avoids certain patterns,
than it avoids also a further pattern. Therefore, it is possible
to apply one of them to a certain class
$S(p_{i_1},p_{i_2},p_{i_3})$ to prove that
$S(p_{i_1},p_{i_2},p_{i_3})\subseteq S(p_1,p_2,p_3,p_4)$ (the
generalization to the case of four and five forbidden pattern is
straightforward). The proof of the first four of them can be found
in \cite{BFP}.

\begin{prop}\label{2}
If $\pi\in S(2-13)$, then $\pi\in S(2-13,21-3)$.
\end{prop}

\begin{prop}\label{3}
If $\pi\in S(31-2)$, then $\pi\in S(31-2,3-12)$.
\end{prop}

\begin{prop}\label{4}
If $\pi\in S(2-31)$, then $\pi\in S(2-31,23-1)$.
\end{prop}

\begin{prop}\label{5}
If $\pi\in S(13-2)$, then $\pi\in S(13-2,1-32)$.
\end{prop}

\begin{prop}\label{u1}
If $\pi\in S(1-23,2-13)$, then $\pi\in S(1-23,2-13,12-3)$.
\end{prop}

\emph{Proof.} Suppose that $\pi$ contains a $12-3$ pattern in the
entries $\pi_i$, $\pi_{i+1}$ and $\pi_k$ ($k>i+1$). Let us
consider the entry $\pi_{i+2}$. It can be neither
$\pi_{i+2}>\pi_{i+1}$ (since $\pi_i\ \pi_{i+1}\ \pi_{i+2}$ would
show a pattern $1-23$) nor $\pi_{i+2}<\pi_{i+1}$ (since
$\pi_{i+1}\ \pi_{i+2}\ \pi_{k}$ would show a pattern $21-3$ which
is forbidden thanks to Proposition \ref{2}).\qed

\noindent(The proof of the following proposition is very similar
and is omitted.)

\begin{prop}\label{u2}
If $\pi\in S(1-23,21-3)$, then $\pi\in S(1-23,21-3,12-3)$.
\end{prop}
%\emph{Proof.} If the entries $\pi_i$, $\pi_{i+1}$ and $\pi_k$ form
%a pattern $12-3$, then it can not be either $\pi_{i+2}>\pi_{i+1}$
%(a pattern $1-23$ would appear in $\pi_i\ \pi_{i+1}\ \pi_{i+2}$)
%or $\pi_{i+2}<\pi_{i+1}$ (a pattern $21-3$ would appear in
%$\pi_{i+1}\ \pi_{i+2}\ \pi_k$).\qed
\begin{prop}\label{u3}
If $\pi\in S(1-23,2-31)$, then $\pi\in S(1-23,2-31,12-3)$.
\end{prop}

\emph{Proof.} Suppose that a pattern $12-3$ appear in $\pi_i$,
$\pi_{i+1}$ and $\pi_k$. If we consider the entry $\pi_{k-1}$,
then it is easily seen that it can be neither
$\pi_i<\pi_{k-1}<\pi_k$ (the entries $\pi_i\ \pi_{k-1}\ \pi_k$
would be $1-23$ pattern like) nor $\pi_{k-1}<\pi_i$ (the entries
$\pi_i\ \pi_{i+1}\ \pi_{k-1}$ would show a pattern $23-1$ which is
forbidden thanks to Proposition \ref{4}). Hence,
$\pi_{k-1}>\pi_k$. We can repeat the same above argument for the
entry $\pi_j$, $j=k-2,k-3,\ldots,i+2$, concluding each time that
$\pi_j>\pi_{j+1}$. When $j=i+2$ a pattern $1-23$ is shown in
$\pi_i\ \pi_{i+1}\ \pi_{i+2}$, which is forbidden.\qed

\begin{prop}\label{u4}
If $\pi\in S(1-23,23-1)$, then $\pi\in S(1-23,23-1,12-3)$.
\end{prop}
This last proposition can be be proved by simply adapting the
argument of the proof of the preceding one.

%\emph{Proof.} Suppose that the entries $\pi_i$, $\pi_{i+1}$ and
%$\pi_k$ form a pattern $12-3$. It must be $\pi_{k-1}>\pi_k$ since
%it can not be either $\pi_i<\pi_{k-1}<\pi_k$ (a pattern $1-23$
%would appear in $\pi_i\ \pi_{k-1}\ \pi_k$) or $\pi_{k-1}<\pi_i$.
%We can repeat the same above argument for the entry $\pi_j$,
%$j=k-2,k-3,\ldots,i+2$, concluding each time that
%$\pi_j>\pi_{j+1}$. When $j=i+2$ a pattern $1-23$ is shown in
%$\pi_i\ \pi_{i+1}\ \pi_{i+2}$, which is forbidden.\qed

\section{Permutations avoiding four patterns}

First of all we recall the results of \cite{BFP} in Tables
\ref{3mot1} and \ref{3mot2}.  For the seek of brevity, for each
symmetry class only a representative is reported. In the first
column of these tables, a name to each symmetry class is given (as
in \cite{BFP}), the second one shows the three forbidden patterns
(the representative) and the third one indicates the sequence
enumerating the permutations avoiding the specified patterns.

Having at our disposal the results for the permutations avoiding
three patterns, the proofs for the case of four forbidden patterns
are conducted following the line indicated in the previous
section. These proofs are all summarized in tables. Tables
\ref{4mot_n}, \ref{4mot_F_n} and \ref{4mot_2^n-1} are related to
the permutations avoiding four patterns enumerated by the
sequences $\{n\}_{n\geq1}$, $\{F_n\}_{n\geq1}$ and
$\{2^{n-1}\}_{n\geq1}$, respectively (the succession $F_n$ denotes
the Fibonacci numbers). As in \cite{BFP}, the empty permutation
with length $n=0$ is not considered, therefore the length is
$n\geq1$. The Tables have to be read as follows: consider the
representative of the symmetry class specified in the rightmost
column of each table; apply the proposition indicated in the
precedent column to the three forbidden patterns which one can
find in Tables \ref{3mot1} and \ref{3mot2} to obtain the four
forbidden patterns written in the column named \emph{avoided
patterns}. At this point, as we explained in the previous section,
the permutations avoiding these four patterns are enumerated by
the same sequence enumerating the permutations avoiding the three
patterns contained in the representative of the symmetry class
indicated in the rightmost column.

The first column of Table \ref{4mot_n} and \ref{4mot_F_n}
specifies a name for the the symmetry class represented by the
four forbidden patterns of the second column. This name is useful
in the next section. Table \ref{4mot_varie} indicates in the first
column the sequence enumerating the permutations avoiding the
patterns of the second column, which are obtained as in the above
tables.

\subsection{Classes enumerated by $\{0\}_{n\geq k}$.}

The classes of four patterns avoiding permutations enumerated by
the sequence $\{0\}_{n\geq k}$ can be handled in a very simple
way. If $S(q_1,q_2,q_3)$, $q_i\in \mathcal M$, is a class of
permutations avoiding three patterns such that
$|S_n(q_1,q_2,q_3)=0|$, for $n\geq k$, then it is easily seen that
$S(q_1,q_2,q_3,r)$, $\forall r\in \mathcal M$, is also enumerated
by the same sequence. Then, each symmetry class from $C1$ to $C7$
(see Table \ref{3mot2}) generates nine symmetry classes by
choosing the pattern $r\neq q_i$, $i=1,2,3$. It is not difficult
to see that all the classes we obtain in this way are not all
different, thanks to the operations of reverse, complement and
reverse-complement. In Table \ref{4mot_0}, only the different
possible cases are presented. Here, the four forbidden patterns
are recovered by adding a pattern of a box of the second column to
the three patterns specified in the box to its right at the same
level (rightmost column). The representative so obtained is
recorded in the leftmost column with a name, which will be useful
in the next section.

\subsection{Classes enumerated by $\{2\}_{n\geq 2}$.}

The enumerating sequences encountered till now (see Tables
\ref{4mot_n}, \ref{4mot_F_n}, \ref{4mot_2^n-1}, \ref{4mot_varie},
\ref{4mot_0}) are all involved in the enumeration of some class of
permutations avoiding three patterns (Tables \ref{3mot1},
\ref{3mot2}). Therefore, applying the eight propositions of the
previous section to the classes of Table \ref{3mot1} and
\ref{3mot2}, the three forbidden patterns have been increased by
one pattern, obtaining Table \ref{4mot_n}, \ref{4mot_F_n},
\ref{4mot_2^n-1}, \ref{4mot_varie} and \ref{4mot_0}. For the
classes enumerated by the sequence $\{2\}_{n\geq 2}$ it is not
possible to use the same strategy, since there are no classes of
permutations avoiding three patterns enumerated by that sequence.
The proofs, in this case, use four easy propositions whose proofs
can be directly derived from the statement of the first four
propositions of the Introduction. We prefer to explicit them the
same.

\begin{prop}\label{contiene_1}
If a permutation $\pi$ contains the pattern $23-1$, then it
contains the pattern $2-31$, too.
\end{prop}

Taking the reverse, the complement and the reverse-complement of
the patterns involved in Prop. \ref{contiene_1}, the following
propositions are obtained:

\begin{prop}\label{contiene_2}
If a permutation $\pi$ contains the pattern $1-32$, then it
contains the pattern $13-2$, too.
\end{prop}

\begin{prop}\label{contiene_3}
If a permutation $\pi$ contains the pattern $21-3$, then it
contains the pattern $2-13$, too.
\end{prop}

\begin{prop}\label{contiene_4}
If a permutation $\pi$ contains the pattern $3-12$, then it
contains the pattern $31-2$, too.
\end{prop}

In Table \ref{4mot_2} the results relating to the enumeration of
the permutations avoiding four patterns enumerated by the sequence
$\{2\}_{n\geq 2}$ (whose proofs are contained in the six next
propositions) are summarized. The four forbidden patterns can be
recovered by choosing one pattern from each column, in the same
box-row of the table.

 In the sequel, $p_i\in A_i$ with $i=1,2,3,4$ where $A_i$ is
a subset of generalized patterns.
\begin{prop}
Let $A_1=\{1-23\}$, $A_2=\{2-31,23-1\}$, $A_3=\{1-32,13-2\}$ and
$A_4=\{3-12,31-2\}$. Then $|S_n(p_1,p_2,p_3,p_4)|=2$ and $S_n=\{n\
(n-1)\ldots 3\ 2\ 1 ,\ (n-1)\ (n-2)\ldots 3\ 2\ 1\ n \}$.
\end{prop}

\emph{Proof.} Let $\sigma\in S_n(p_2,p_3)$. Then, $\sigma_1=n$ or
$\sigma_n=n$, otherwise, if $\sigma_i=n$ with $i\neq 1,n\ $, the
entries $\sigma_{i-1}\sigma_{i}\sigma_{i+1}$ would be a forbidden
pattern $p_2$ or $p_3$.

If $\rho\in S_n(p_1,p_3)$, then $\rho_{n-1}=1$ or $\rho_n=1$,
otherwise, if $\rho_i=1$ with $i<n-1$, then the entries
$\rho_i\rho_{i+1}\rho_{i+2}$, would be a forbidden pattern $p_1$
or $p_3$.

Therefore, if $\pi\in S_n(p_1,p_2,p_3)$, then there are only the
following three cases for $\pi$:
\begin{enumerate}
    \item $\pi_n=n$ and $\pi_{n-1}=1$. In this case $\pi=(n-1)\ (n-2) \ldots
    2\ 1\ n$, otherwise, if an ascent appears in $\pi_j\pi_{j+1}$ with
    $j=1,2,\ldots, n-3$, the entries $\pi_j\pi_{j+1}\pi_{n-1}$
    would show the pattern $23-1$ and $\pi$ would contain the pattern
    $2-31$, too (see Prop. \ref{contiene_1}).
    \item $\pi_1=n$ and $\pi_n=1$. In this case $\pi=n\ (n-1)\ldots 3\ 2\
    1$, otherwise, if an ascent appears in $\pi_j\pi_{j+1}$ with
    $j=2,3,\ldots, n-2$, the entries $\pi_j\pi_{j+1}\pi_n$ would
    show the pattern $23-1$ and $\pi$ would contain the pattern
    $2-31$, too (see Prop. \ref{contiene_1}).
    \item $\pi_1=n$ and $\pi_{n-1}=1$ (and $\pi_n=k<n$).
\end{enumerate}

If $\pi$ has to avoid the pattern $p_4$, too ($\pi\in
S_n(p_1,p_2,p_3,p_4)$), then the third above case is not allowed
since $\pi_1\pi_{n-1}\pi_n$ are a $3-12$ pattern which induces an
occurrence of $31-2$ in $\pi$ (Prop. \ref{contiene_4}). \qed

\begin{prop}
Let $A_1=\{1-23\}$, $A_2=\{2-13,21-3\}$, $A_3=\{1-32,13-2\}$ and
$A_4=\{3-12,31-2\}$. Then $|S_n(p_1,p_2,p_3,p_4)|=2$ and $S_n=\{n\
(n-1)\ldots 3\ 2\ 1,\ (n-1)\ n\ (n-2)\ (n-3)\ldots 2\ 1\}$.
\end{prop}

\emph{Proof.} If $\sigma\in S_n(p_1,p_2)$, then $\pi_1=n$ or
$\pi_2=n$. If $\rho\in S_n(p_1,p_3)$, then $\pi_n=1$ or
$\pi_{n-1}=1$. Then, if $\pi\in S_n(p_1,p_2,p_3)$, there are only
the four following cases:
\begin{enumerate}
    \item $\pi_1=n$ and $\pi_n=1$.
    \item $\pi_2=n$ and $\pi_n=1$. In this case $\pi_1=n-1$,
    otherwise if $\pi_k=n-1$ with $k>3$, then $\pi_{k-2}\pi_{k-1}\pi_k$ is
    a $1-23$ pattern or a $21-3$ pattern which induces an occurrence of
    $2-13$ (Prop. \ref{contiene_3}). If $k=3$, then $\pi_1\pi_2\pi_3$ is a
    $1-32$ or $13-2$ pattern which are forbidden.
    \item $\pi_1=n$ and $\pi_{n-1}=1$.
    \item $\pi_2=n$ and $\pi_{n-1}=1$. For the same reasons of
    case 2, it is $\pi_1=n-1$.
\end{enumerate}

If $\pi$ has to avoid $p_4$, too ($\pi\in S_n(p_1,p_2,p_3,p_4)$),
then the third and the fourth above cases are not allowed since
$\pi_1\pi_{n-1}\pi_n$ are a $3-12$ pattern which induces an
occurrence of $31-2$ (Prop. \ref{contiene_4}). Moreover, the
permutations of the above cases $1$ and $2$, must be such that
there are not ascents $\pi_i\pi_{i+1}$ between $n$ and $1$ in
order to avoid $p_4$. Then, $\pi=n\ (n-1)\ldots 3\ 2\ 1 $ or
$\pi=(n-1)\ n\ (n-2)\ldots 3\ 2\ 1$. \qed

\begin{prop}\label{16classi}
Let $A_1=\{2-13,21-3\}$, $A_2=\{2-31,23-1\}$, $A_3=\{1-32,13-2\}$
and $A_4=\{3-12,31-2\}$. Then $|S_n(p_1,p_2,p_3,p_4)|=2$ and
$S_n=\{n\ (n-1)\ldots 2\ 1,\ 1\ 2\ldots n\}$.
\end{prop}

\emph{Proof.} It is easily seen that each three consecutive
elements of $\pi$ can only be in increasing or decreasing order.
\qed

\begin{prop}
Let $A_1=\{12-3\}$, $A_2=\{2-13,21-3\}$, $A_3=\{2-31,23-1\}$ and
$A_4=\{32-1\}$. Then $|S_n(p_1,p_2,p_3,p_4)|=2$ and $S_n=\{1\ n\
2\ (n-1)\ldots,\ n\ 1\ (n-1)\ 2\ldots\}$.
\end{prop}

\emph{Proof.} If $\pi\in S_n(p_1,p_2,p_3,p_4)$, then it is easy to
see that $\pi_1\pi_2=1\ n$ or $\pi_1\pi_2=n\ 1$. Considering the
sub-permutation $\pi_2\pi_3\ldots \pi_n$, in the same way we
deduce $\pi_2\pi_3=2\ (n-1)$ or $\pi_2\pi_3=(n-1)\ 2$. The thesis
follows by recursively using the above argument. \qed

\begin{prop}
Let $A_1=\{1-23\}$, $A_2=\{2-13,21-3\}$, $A_3=\{2-31,23-1\}$ and
$A_4=\{3-12,31-2\}$. Then $|S_n(p_1,p_2,p_3,p_4)|=2$ and $S_n=\{n\
(n-1) \ldots 1,\ 1\ n\ (n-1) \ldots 3\ 2\}$.
\end{prop}

\emph{Proof.} Let $\pi\in S_n(p_1,p_2,p_3,p_4)$. It is $\pi_1=n$
or $\pi_2=n$, otherwise a $1-23$ or $p_2$ pattern would appear.

If $\pi_1=n$, then $\pi=n\ (n-1) \ldots 1$ since if an ascent
appears in $\pi_i\pi_{i+1}$, the entries $\pi_1\pi_i\pi_{i+1}$ are
a $p_4$ pattern.

If $\pi_2=n$, then $\pi_1=1$ since the $p_3$ pattern has to be
avoided. Moreover, in this case, it is $\pi_j>\pi_{j+1}$ with
$j=3,4,\ldots,(n-1)$ in order to avoid $1-23$. Then $\pi=1\ n\
(n-1)\ldots 2\ 1$. \qed

\begin{prop}
Let $A_1=\{1-23\}$, $A_2=\{2-13,21-3\}$, $A_3=\{2-31,23-1\}$ and
$A_4=\{1-32,13-2\}$. Then $|S_n(p_1,p_2,p_3,p_4)|=2$ and $S_n=\{n\
(n-1) \ldots 3\ 2\ 1,\ n\ (n-1) \ldots 3\ 1\ 2\}$.
\end{prop}

\emph{Proof.} Let $\pi\in S_n(p_1,p_2,p_3,p_4)$. The entries $1$
and $2$ have to be adjacent in order to avoid $p_3$ and $p_4$ and
$\pi_n=1$ or $\pi_{n-1}=1$ in order to avoid $p_1$ and $p_4$. So,
$\pi_{n-1}\pi_n=1\ 2$ or $\pi_{n-1}\pi_n=2\ 1$. Moreover, each
couple of adjacent elements $\pi_j\pi_{j+1}$ must be a descent,
otherwise a $23-1$ pattern (which induces an occurrences of
$2-31$) would appear. Then $\pi=n\ (n-1) \ldots 3\ 2\ 1$ or $\pi=
n\ (n-1) \ldots 3\ 1\ 2$. \qed

The conjecture stated in \cite{CM} about the permutations
enumerated by $\{2\}_{n\geq2}$ declares that there are $42$
symmetry classes of such permutations, while from Table
\ref{4mot_2} it is possible to deduce $52$ symmetry classes.
Nevertheless, it is not difficult to check that these classes are
not all different: for example the symmetry class
$\{2-13,2-31,1-32,31-2\}$ is the same of $\{2-13,23-1,13-2,31-2\}$
(the second one is the reverse of the first one). Note that the
repetitions come out only from the third box-row of Table
\ref{4mot_2}.

\section{Permutations avoiding five patterns}

\subsection{Classes enumerated by $\{1\}_{n\geq 1}$}

The sequence $\{1\}_{n\geq 1}$ does not enumerate any class of
permutations avoiding four patterns, so that we can not apply the
same method of the previous section using the proposition of the
Introduction.

Referring to Proposition \ref{16classi}, we deduce that there are
sixteen different classes $S_n(p_1,p_2,p_3,p_4)$ such that $p_i
\in A_i$ with $i=1,2,3,4$. We recall that
$|S_n(p_1,p_2,p_3,p_4)|=2$ and $S_n(p_1,p_2,p_3,p_4)=\{n\
(n-1)\ldots 2\ 1,\ 1\ 2\ldots n\}$. If a permutation $\pi\in
S_n(p_1,p_2,p_3,p_4)$ has to avoid the pattern $1-23$, too, then
$\pi=n\ (n-1)\ldots 2\ 1$ and $|S_n(p_1,p_2,p_3,p_4,1-23)|=1$.

Then, it is easy to see that the five forbidden patterns avoided
by the permutations enumerated by $\{1\}_{n\geq 1}$ can be
recovered by considering the four patterns chosen from the third
box-row of Table \ref{4mot_2} (one pattern from each column) and
the pattern $1-23$. We do not present the relative table.

\subsection{Classes enumerated by $\{0\}_{n\geq k}$}

This case is treated as the case of the permutations avoiding four
patterns. It is sufficient to add a pattern $r\in\mathcal M$ to
each representative (from O1 to O37 in Table \ref{4mot_0}) of four
forbidden patterns of Table \ref{4mot_0} in order to obtain a
representative $T$ of five forbidden patterns such that
$|S_n(T)|=0, n\geq4$. In Table \ref{5mot_0} we present the
different representatives $T$ which can be derived from Table
\ref{4mot_0}. The five forbidden patterns of each representative
are a pattern chosen in a box of the first column and the four
patterns indicated by the representative (which refer to Table
\ref{4mot_0}) in the second box at the same level. In the table,
only the different representatives of five patterns are presented.

\subsection{Classes enumerated by $\{2\}_{n\geq k}$,
$\{n\}_{n\geq 1}$, $\{F_n\}_{n\geq 1}$}

Tables \ref{5mot_2(1)} and \ref{5mot_2(2)} summarize the results
related to the permutations avoiding five patterns enumerated by
$\{2\}_{n\geq k}$. The five forbidden patterns are obtained by
considering a representative of four forbidden patterns of the
rightmost column and the pattern specified in the corresponding
box of the preceding column. The first column indicates which is
the proposition to apply. Note that each representative of four
patterns (rightmost column) can be found in Table \ref{4mot_2}.

The reading of Tables \ref{5mot_n} and \ref{5mot_F_n} (related to
the sequences $\{n\}_{n\geq 1}$ and $\{F_n\}_{n\geq 1}$,
respectively) is as usual: apply the proposition specified in the
first column to recover the representative of five forbidden
patterns which is composed by the pattern of the second column and
the four patterns of the representative indicated in the rightmost
column. Here, the names of the representatives refer to Tables
\ref{4mot_n} and \ref{4mot_F_n}.

%%%%%%%%%%
%%%%%%%%%%
%%%%%%%%%%
%%SEZIONE[TAVOLE]
%%%%%%%%%%
%%%%%%%%%%

%PRIMA TAVOLA DEI 3 MOTIVI
\section{Tables}\label{tables}
\begin{table}[!h]
\begin{center}
\begin{tabular}{|c|c|c|}
symmetry class & avoided patterns & enumerating sequence\\
\hline \hline

N1&\{1-23,2-13,3-12\}&\\
%\hline
N2&\{1-23,2-13,31-2\}&\\
%\hline
N3&\{1-23,21-3,3-12\}&\\
%\hline
N4&\{1-23,21-3,31-2\}&\\
%\hline
N5&\{12-3,3-12,2-13\}&\\
%\hline
N6&\{12-3,3-12,21-3\}&\\
%\hline
N7&\{12-3,31-2,2-13\}&\\
%\hline
N8&\{12-3,31-2,21-3\}&\\
%\hline
N9&\{1-23,2-13,2-31\}&\\
%\hline
N10&\{1-23,2-13,23-1\}&\\
%\hline
N11&\{1-23,21-3,2-31\}&\\
%\hline
N12&\{1-23,21-3,23-1\}&\\
%\hline
N13&\{2-13,2-31,1-32\}& $\{n\}_{n\geq 1}$\\
%\hline
N14&\{2-13,23-1,1-32\}&\\
%\hline
N15&\{2-13,2-31,13-2\}&\\
%\hline
N16&\{2-13,23-1,13-2\}&\\
%\hline
N17&\{2-31,21-3,13-2\}&\\
%\hline
N18&\{2-31,21-3,1-32\}&\\
%\hline
N19&\{13-2,21-3,23-1\}&\\
%\hline
N20&\{21-3,23-1,1-32\}&\\
%\hline
N21&\{1-23,2-31,31-2\}&\\
%\hline
N22&\{1-23,23-1,31-2\}&\\
%\hline
N23&\{1-23,2-31,3-12\}&\\
%\hline
N24&\{1-23,1-32,3-21\}&\\
\hline
A1&\{1-23,12-3,23-1\}&\\
A2&\{2-31,23-1,1-32\}&\\
A3&\{2-31,23-1,13-2\}&\\
A4&\{1-23,12-3,2-13\}&\\
A5&\{1-23,2-13,21-3\}&$\{2^{n-1}\}_{n\geq 1}$\\
A6&\{1-23,3-12,31-2\}&\\
A7&\{31-2,3-12,13-2\}&\\
A8&\{31-2,3-12,1-32\}&\\
A9&\{2-13,21-3,1-32\}&\\
A10&\{2-13,21-3,13-2\}&\\
\hline
A11&\{1-23,23-1,3-12\}&$\{2{n-2}+1\}_{n\geq 1}$\\
\hline
\end{tabular}
\end{center}
\caption{permutations avoiding three patterns}\label{3mot1}
\end{table}

\begin{table}
\begin{center}
\begin{tabular}{|c|c|c|}
symmetry class & avoided patterns & enumerating sequence\\
\hline \hline

F1&$\{1-23,2-13,1-32\}$&\\
F2&$\{1-23,2-13,13-2\}$&\\
F3&$\{1-23,21-3,13-2\}$&\\
F4&$\{1-23,13-2,3-12\}$& $\{F_n\}_{n\geq1}$\\
F5&$\{1-23,1-32,3-12\}$&\\
F6&$\{1-23,1-32,31-2\}$&\\
F7&$\{1-23,13-2,31-2\}$&\\
\hline
M1&$\{1-23,12-3,21-3\}$&\\
M2&$\{12-3,21-3,2-13\}$& $\{M_n\}_{n\geq1}$\\
\hline
B1&$\{1-23,21-3,1-32\}$& $\{{n \choose \lceil n/2 \rceil}\}_{n\geq1}$\\
\hline
B2&$\{12-3,1-23,31-2\}$&\\
B3&$\{1-23,2-31,23-1\}$& $\{1+{n \choose 2}\}$\\
\hline
C8&$\{12-3,2-13,32-1\}$& $\{3\}_{n\geq 3}$\\
\hline
C1&$\{1-23,2-13,3-21\}$&\\
C2&$\{1-23,23-1,32-1\}$&\\
C3&$\{1-23,2-13,32-1\}$&\\
C4&$\{1-23,12-3,3-21\}$& $\{0\}_{n\geq k}$\\
C5&$\{1-23,21-3,3-21\}$&\\
C6&$\{1-23,21-3,32-1\}$&\\
C7&$\{1-23,2-31,32-1\}$&\\
\hline
\end{tabular}
\end{center}
\caption{permutations avoiding three patterns}\label{3mot2}
\end{table}

\begin{table}
\begin{center}
\begin{tabular}{|c|c|c|c|}
\hline

\multicolumn{4}{|c|}{\textbf{Enumerating sequence: $\{n\}_{n\geq 1}$}}\\

\hline \hline

\emph{name} & \emph{avoided patterns} & \emph{apply Proposition} & \emph{to the symmetry class}\\

\hline $d1$ & $\{1-23,2-13,3-12,21-3\}$ & \ref{2}  & N1 \\
\hline $d2$ & $\{1-23,2-13,31-2,21-3\}$ & \ref{2}  & N2 \\
\hline $d3$ & $\{1-23,2-13,31-2,3-12\}$ & \ref{3}  & N2 \\
\hline $d4$ & $\{1-23,21-3,31-2,3-12\}$ & \ref{3}  & N4 \\
\hline $d5$ & $\{12-3,3-12,2-13,21-3\}$ & \ref{2}  & N5 \\
\hline $d6$ & $\{12-3,31-2,2-13,21-3\}$ & \ref{2}  & N7 \\
\hline $d7$ & $\{12-3,31-2,2-13,3-12\}$ & \ref{3}  & N7 \\
\hline $d8$ & $\{12-3,31-2,21-3,3-12\}$ & \ref{3}  & N8 \\
\hline $d9$ & $\{1-23,2-13,2-31,21-3\}$ & \ref{2}  & N9 \\
\hline $d10$ & $\{1-23,2-13,2-31,23-1\}$ & \ref{4}  & N9 \\
\hline $d11$ & $\{1-23,2-13,23-1,21-3\}$ & \ref{2}  & N10\\
\hline $d12$ & $\{1-23,21-3,2-31,23-1\}$ & \ref{4}  & N11\\
\hline $d13$ & $\{2-13,2-31,1-32,21-3\}$ & \ref{2}  & N13\\
\hline $d14$ & $\{2-13,2-31,1-32,23-1\}$ & \ref{4}  & N13\\
\hline $d15$ & $\{2-13,23-1,1-32,21-3\}$ & \ref{2}  & N14\\
\hline $d16$ & $\{2-13,2-31,13-2,21-3\}$ & \ref{2}  & N15\\
\hline $d17$ & $\{2-13,2-31,13-2,23-1\}$ & \ref{4}  & N15\\
\hline $d18$ & $\{2-13,2-31,13-2,1-32\}$ & \ref{5}  & N15\\
\hline $d19$ & $\{2-13,23-1,13-2,21-3\}$ & \ref{2}  & N16\\
\hline $d20$ & $\{2-13,23-1,13-2,1-32\}$ & \ref{5}  & N16\\
\hline $d21$ & $\{2-31,21-3,13-2,23-1\}$ & \ref{4}  & N17\\
\hline $d22$ & $\{2-31,21-3,13-2,1-32\}$ & \ref{5}  & N17\\
\hline $d23$ & $\{2-31,21-3,1-32,23-1\}$ & \ref{4}  & N18\\
\hline $d24$ & $\{13-2,21-3,23-1,1-32\}$ & \ref{5}  & N19\\
\hline $d25$ & $\{1-23,2-31,31-2,23-1\}$ & \ref{4}  & N21\\
\hline $d26$ & $\{1-23,2-31,31-2,3-12\}$ & \ref{3}  & N21\\
\hline $d27$ & $\{1-23,23-1,31-2,3-12\}$ & \ref{3}  & N22\\
\hline $d28$ & $\{1-23,2-31,3-12,23-1\}$ & \ref{4}  & N23\\
\hline $d29$ & $\{1-23,2-13,31-2,12-3\}$ & \ref{u1} & N2 \\
\hline $d30$ & $\{1-23,2-13,3-12,12-3\}$ & \ref{u1} & N1 \\
\hline $d31$ & $\{1-23,2-13,2-31,12-3\}$ & \ref{u1} & N9 \\
\hline $d32$ & $\{1-23,2-13,23-1,12-3\}$ & \ref{u1} & N10\\
\hline $d33$ & $\{1-23,21-3,2-31,12-3\}$ & \ref{u2} & N11\\
\hline $d34$ & $\{1-23,21-3,23-1,12-3\}$ & \ref{u2} & N12\\
\hline $d35$ & $\{1-23,21-3,31-2,12-3\}$ & \ref{u2} & N4 \\
\hline $d36$ & $\{1-23,21-3,3-12,12-3\}$ & \ref{u2} & N3 \\
\hline $d37$ & $\{1-23,2-31,3-12,12-3\}$ & \ref{u3} & N23\\
\hline $d38$ & $\{1-23,2-31,31-2,12-3\}$ & \ref{u3} & N21\\
\hline
\end{tabular}
\end{center}
\caption{permutations avoiding four patterns}\label{4mot_n}
\end{table}

\begin{table}
\begin{center}
\begin{tabular}{|c|c|c|c|}
\hline

\multicolumn{4}{|c|}{\textbf{Enumerating sequence: $\{F_n\}_{n\geq 1}$}}\\

\hline \hline

\emph{name} & \emph{avoided patterns} & \emph{apply Proposition} & \emph{to the symmetry class}\\

\hline $e1$ & $\{1-23,2-13,1-32,21-3\}$ & \ref{2}  & F1\\
\hline $e2$ & $\{1-23,2-13,1-32,12-3\}$ & \ref{u1} & F1\\
\hline $e3$ & $\{1-23,2-13,13-2,21-3\}$ & \ref{2}  & F2\\
\hline $e4$ & $\{1-23,2-13,13-2,1-32\}$ & \ref{5}  & F2\\
\hline $e5$ & $\{1-23,2-13,13-2,12-3\}$ & \ref{u1} & F2\\
\hline $e6$ & $\{1-23,21-3,13-2,1-32\}$ & \ref{5}  & F3\\
\hline $e7$ & $\{1-23,13-2,3-12,1-32\}$ & \ref{5}  & F4\\
\hline $e8$ & $\{1-23,1-32,31-2,3-12\}$ & \ref{3}  & F6\\
\hline $e9$ & $\{1-23,13-2,31-2,1-32\}$ & \ref{5}  & F7\\
\hline $e10$ & $\{1-23,13-2,31-2,3-12\}$ & \ref{3}  & F7\\
\hline
\end{tabular}
\end{center}
\caption{permutations avoiding four patterns}\label{4mot_F_n}
\end{table}

\begin{table}
\begin{center}
\begin{tabular}{|c|c|c|}
\hline

\multicolumn{3}{|c|}{\textbf{Enumerating sequence: $\{2^{n-1}\}_{n\geq 1}$}}\\

\hline \hline

\emph{avoided patterns} & \emph{apply Proposition} & \emph{to the symmetry class}\\

\hline $\{1-23,12-3,2-13,21-3\}$ & \ref{2}  & A4 \\
\hline $\{31-2,3-12,13-2,1-32\}$ & \ref{5}  & A7 \\
\hline $\{2-13,21-3,13-2,1-32\}$ & \ref{5}  & A10\\
\hline $\{2-31,23-1,1-32,13-2\}$ & \ref{5}  & A3 \\
\hline
\end{tabular}
\end{center}
\caption{permutations avoiding four patterns}\label{4mot_2^n-1}
\end{table}

\begin{table}
\begin{center}
\small{
\begin{tabular}{|c|c|c|c|}
\hline
$
\begin{array}{c}
  \textbf{Enumerating} \\
  \textbf{sequence} \\
\end{array}$
& \emph{avoided patterns} &
\emph{apply Proposition} & \emph{to the symmetry class} \\
\hline \hline
$\{1+{n \choose 2}\}_{n\geq 1}$ & $\{12-3,1-23,31-2,3-12\}$ & \ref{3}  & B2 \\
\hline
$\{{n \choose \lceil n/2 \rceil}\}_{n\geq 1}$ & $\{1-23,21-3,1-32,12-3\}$ & \ref{u2} & B1 \\
\hline
$\{2^{n-2}+1\}_{n\geq 1}$ & $\{1-23,23-1,3-12,12-3\}$ & \ref{u4} & A11 \\
\hline
$\{3\}_{n\geq 3}$ & $\{12-3,2-13,32-1,21-3\}$ & \ref{2} & C8 \\
\hline
\end{tabular}
}
\end{center}
\caption{permutations avoiding four patterns}\label{4mot_varie}
\end{table}

\begin{table}
\begin{center}
\begin{tabular}{|c|c|c|}
\hline
\multicolumn{3}{|c|}{\textbf{Enumerating sequence: $\{0\}_{n\geq k}$}}\\
\hline\hline

\emph{name}&\emph{choose a pattern from the following to add} &\emph{to the symmetry class} \\
\hline \hline

O1& $12-3$&\\
O2& $1-32$&\\
O3& $13-2$&\\
O4& $3-12$&\\
O5& $31-2$&$\{1-23,2-13,3-21\}$ (C1)\\
O6& $21-3$&\\
O7& $2-31$&\\
O8& $23-1$&\\
O9& $32-1$&\\
\hline
O10& $12-3$&\\
O11& $1-32$&\\
O12& $13-2$&\\
O13& $3-12$&\\
O14& $31-2$&$\{1-23,23-1,32-1\}$ (C2)\\
O15& $2-13$&\\
O16& $21-3$&\\
O17& $2-31$&\\
O18& $3-21$&\\
\hline
O19& $12-3$&\\
O20& $13-2$&\\
O21& $3-12$&\\
O22& $31-2$&$\{1-23,2-13,32-1\}$ (C3)\\
O23& $21-3$&\\
O24& $2-31$&\\
\hline
O25& $31-2$&\\
O26& $1-32$&\\
O27& $23-1$&$\{1-23,12-3,3-21\}$ (C4)\\
O28& $32-1$&\\
\hline
O29& $1-32$&\\
O30& $13-2$&\\
O31& $3-12$&$\{1-23,21-3,3-21\}$ (C5)\\
O32& $31-2$&\\
O33& $23-1$&\\
\hline
O34& $13-2$&\\
O35& $3-12$&$\{1-23,21-3,32-1\}$ (C6)\\
O36& $2-31$&\\
\hline
O37& $13-2$&$\{1-23,2-31,32-1\}$ (C7)\\
\hline
\end{tabular}
\end{center}
\caption{permutations avoiding four patterns}\label{4mot_0}
\end{table}

%TABELLA SEQUENZA 2

\begin{table}
\begin{center}
\begin{tabular}{|c|c|c|c|}
\hline
\multicolumn{4}{|c|}{\textbf{Enumerating sequence: $\{2\}_{n\geq 2}$}}\\
\hline \hline

\emph{1st pattern}&\emph{2nd pattern}&\emph{3rd pattern}&\emph{4th
pattern}\\
\hline \hline

$\begin{array}{c}
  1-23 \\
\end{array}$

&

$\begin{array}{l}
  2-31\ \ or \\
  23-1 \\
\end{array}$

&

$\begin{array}{l}
  1-32\ \ or \\
  13-2 \\
\end{array}$

&

$\begin{array}{l}
  3-12\ \ or \\
  31-2 \\
\end{array}$

\\ \hline

$\begin{array}{c}
  1-23 \\
\end{array}$

&

$\begin{array}{l}
  2-13\ \ or \\
  21-3 \\
\end{array}$

&

$\begin{array}{l}
  1-32\ \ or \\
  13-2 \\
\end{array}$

&

$\begin{array}{l}
  3-12\ \ or \\
  31-2 \\
\end{array}$

\\ \hline

$\begin{array}{l}
  2-13\ \ or \\
  21-3 \\
\end{array}$

&

$\begin{array}{l}
  2-31\ \ or \\
  23-1 \\
\end{array}$

&

$\begin{array}{l}
  1-32\ \ or \\
  13-2 \\
\end{array}$

&

$\begin{array}{l}
  3-12\ \ or \\
  31-2 \\
\end{array}$

\\ \hline

$\begin{array}{c}
  12-3 \\
\end{array}$

&

$\begin{array}{l}
  2-13\ \ or \\
  21-3 \\
\end{array}$

&

$\begin{array}{l}
  2-31\ \ or \\
  23-1 \\
\end{array}$

&

$\begin{array}{c}
  32-1 \\
\end{array}$

\\ \hline

$\begin{array}{c}
  1-23 \\
\end{array}$

&

$\begin{array}{l}
  2-13\ \ or \\
  21-3 \\
\end{array}$

&

$\begin{array}{l}
  2-31\ \ or \\
  23-1 \\
\end{array}$

&

$\begin{array}{l}
  3-12\ \ or \\
  31-2 \\
\end{array}$

\\ \hline

$\begin{array}{c}
  1-23 \\
\end{array}$

&

$\begin{array}{l}
  2-13\ \ or \\
  21-3 \\
\end{array}$

&

$\begin{array}{l}
  2-31\ \ or \\
  23-1\\
\end{array}$

&

$\begin{array}{l}
  1-32\ \ or \\
  13-2 \\
\end{array}$

\\ \hline
\end{tabular}
\end{center}
\caption{permutations avoiding four patterns}\label{4mot_2}
\end{table}

\begin{table}[!h]
\begin{center}
\begin{tabular}{|l|c|}
\hline
\multicolumn{2}{|c|}{\textbf{Enumerating sequence: $\{0\}_{n\geq k}$}}\\
\hline\hline

\emph{choose a pattern from the following to add} &\emph{to the symmetry class} \\
\hline

 $\begin{array}{l}
  21-3,\ 2-31,\ 23-1,\ 1-32,\ 13-2,\\
  3-12,\ 31-2,\ 32-1\\
\end{array}$
&O1\\
\hline

 $\begin{array}{l}
  2-31,\ 23-1,\ 1-32,\ 13-2,\ 3-12,\\
  31-2\\
\end{array}$
&O6\\
\hline

 $\begin{array}{l}
  21-3,\ 2-31,\ 23-1,\ 1-32,\ 3-12,\\
  31-2\\
\end{array}$
&O19\\
\hline

 $\begin{array}{l}
  2-31,\  23-1,\  1-32,\  13-2,\ 3-12 \\
  31-2\\
\end{array}$
&O23\\
\hline

 $\begin{array}{l}
  1-32, 13-2, 3-12, 31-2,\ 32-1 \\
  \end{array}$
&O8\\
\hline

 $\begin{array}{l}
  23-1,\  1-32,\  13-2,\  31-2,\  3-12\\
\end{array}$
&O24\\
\hline

 $\begin{array}{l}
  12-3,\  32-1,\  13-2,\  3-12,\  31-2 \\

\end{array}$
&O29\\
\hline

 $\begin{array}{l}
  3-12,\  13-2,\  1-32,\  23-1,\  12-3\\
\end{array}$
&O36\\
\hline

 $\begin{array}{l}
  1-32,\  13-2,\  3-12,\  31-2\\
\end{array}$
&O15\\
\hline

 $\begin{array}{l}
  13-2,\ 3-12,\ 31-2\\
\end{array}$
&O2\\
\hline

 $\begin{array}{l}
  1-32,\ 2-31,\ 31-2\\
\end{array}$
&O10\\
\hline

 $\begin{array}{l}
  12-3,\ 13-2,\ 3-12\\
\end{array}$
&O32\\
\hline

 $\begin{array}{l}
  1-32,\ 13-2,\ 32-1\\
\end{array}$
&O33\\
\hline

 $\begin{array}{l}
  3-21,\ 23-1,\ 1-32\\
\end{array}$
&O34\\
\hline

 $\begin{array}{l}
 3-12,\ 31-2\\
\end{array}$
&O3\\
\hline

 $\begin{array}{l}
  1-32,\ 13-2\\
\end{array}$
&O7\\
\hline

 $\begin{array}{l}
  13-2,\ 3-12\\
\end{array}$
&O9\\
\hline

 $\begin{array}{l}
  21-3,\ 13-2\\
\end{array}$
&O11\\
\hline

 $\begin{array}{l}
  1-32,\ 3-12\\
\end{array}$
&O20\\
\hline

 $\begin{array}{l}
  3-12,\ 23-1\\
\end{array}$
&O26\\
\hline

 $\begin{array}{l}
  2-31,\ 32-1\\
\end{array}$
&O27\\
\hline

 $\begin{array}{l}
 3-12\\
\end{array}$
&O5\\
\hline

 $\begin{array}{l}
  2-31\\
\end{array}$
&O12\\
\hline

 $\begin{array}{l}
  1-32\\
\end{array}$
&O21\\
\hline

 $\begin{array}{l}
  3-12\\
\end{array}$
&O30\\
\hline

 $\begin{array}{l}
  23-1\\
\end{array}$
&O35\\
\hline
\end{tabular}
\end{center}
\caption{permutations avoiding five patterns}\label{5mot_0}
\end{table}

\begin{table}
\begin{center}
\begin{tabular}{|c|c|l|}
\hline
\multicolumn{3}{|c|}{\textbf{Enumerating sequence: $\{2\}_{n\geq k}$}}\\
\hline\hline

\emph{thanks to Proposition} &\emph{add the pattern} &\emph{to the patterns} \\
\hline

\ref{2} & $21-3$ & $\begin{array}{l}
  \{1-23, 2-13, 1-32, 3-12\}\ or \\
  \{1-23, 2-13, 1-32, 31-2\}\ or \\
  \{1-23, 2-13, 13-2, 3-12\}\ or \\
  \{1-23, 2-13, 1-32, 31-2\} \\
\end{array}$\\
\hline

\ref{2} & $21-3$ & $\begin{array}{l}
  \{1-23, 2-13, 2-31, 1-32\}\ or \\
  \{1-23 , 2-13 , 2-31 , 13-2\}\ or \\
  \{1-23 , 2-13 , 23-1 , 1-32\}\ or \\
  \{1-23 , 2-13 , 23-1 , 13-2\} \\
\end{array}$\\
\hline

\ref{2} & $21-3$ & $\begin{array}{l}
  \{1-23 , 2-13 , 2-31 , 3-12\}\ or \\
  \{1-23 , 2-13 , 2-31 , 31-2\}\ or \\
  \{1-23 , 2-13 , 23-1 , 3-12\}\ or \\
  \{1-23 , 2-13 , 23-1 , 31-2\} \\
\end{array}$\\
\hline

\ref{2} & $21-3$ & $\begin{array}{l}
  \{1-23 , 2-13 , 23-1 , 32-1\}\ or \\
  \{1-23 , 2-13 , 2-31 , 32-1\} \\
\end{array}$\\
\hline

\ref{2} & $21-3$ & $\begin{array}{l}
  \{2-13 , 2-31 , 1-32 , 3-12\}\ or \\
  \{2-13 , 2-31 , 1-32 , 31-2\}\ or \\
  \{2-13 , 2-31 , 13-2 , 31-2\} \\
\end{array}$\\
\hline

\ref{2} & $21-3$ & $\begin{array}{l}
  \{2-13 , 23-1 , 1-32 , 31-2\}\ or \\
  \{2-13 , 23-1 , 1-32 , 3-12\}\ or \\
  \{2-13 , 23-1 , 13-2 , 31-2\}\ or \\
  \{2-13 , 23-1 , 13-2 , 3-12\} \\
\end{array}$\\
\hline

\ref{3} & $3-12$ & $\begin{array}{l}
  \{1-23 , 2-13 , 2-31 , 31-2\}\ or \\
  \{1-23 , 2-13 , 23-1 , 31-2\} \\
\end{array}$\\
\hline

\ref{3} & $3-12$ & $\begin{array}{l}
  \{1-23 , 2-13 , 1-32 , 31-2\}\ or \\
  \{1-23 , 2-13 , 13-2 , 31-2\} \\
\end{array}$\\
\hline

\ref{3} & $3-12$ & $\begin{array}{l}
  \{1-23 , 2-13 , 2-31 , 31-2\}\ or \\
  \{1-23 , 2-13 , 23-1 , 31-2\} \\
\end{array}$\\
\hline

\ref{3} & $3-12$ & $\begin{array}{l}
  \{1-23 , 21-3 , 1-32 , 31-2\}\ or \\
  \{1-23 , 21-3 , 13-2 , 31-2\} \\
\end{array}$\\
\hline

\ref{3} & $3-12$ & $\begin{array}{l}
  \{1-23 , 2-31 , 1-32 , 31-2\}\ or \\
  \{1-23 , 2-31 , 13-2 , 31-2\} \\
\end{array}$\\
\hline

\ref{3} & $3-12$ & $\begin{array}{l}
  \{1-23 , 23-1 , 1-32 , 31-2\}\ or \\
  \{1-23 , 23-1 , 13-2 , 31-2\} \\
\end{array}$\\
\hline

\ref{4} & $23-1$ & $\begin{array}{l}
  \{2-13 , 2-31 , 1-32 , 31-2\} \\
\end{array}$\\
\hline

\ref{4} & $23-1$ & $\begin{array}{l}
  \{1-23 , 2-31 , 13-2 , 3-12\}\ or \\
  \{1-23 , 2-31 , 13-2 , 31-2\} \\
\end{array}$\\
\hline

\end{tabular}
\end{center}
\caption{permutations avoiding five patterns}\label{5mot_2(1)}
\end{table}

\begin{table}
\begin{center}
\begin{tabular}{|c|c|l|}
\hline
\multicolumn{3}{|c|}{\textbf{Enumerating sequence: $\{2\}_{n\geq k}$}}\\
\hline\hline

\emph{thanks to Proposition} &\emph{add the pattern} &\emph{to the patterns} \\
\hline

\ref{4} & $23-1$ & $\begin{array}{l}
  \{1-23 , 2-31 , 1-32 , 3-12\}\ or \\
  \{1-23 , 2-31 , 1-32 , 31-2\} \\
\end{array}$\\
\hline

\ref{4} & $23-1$ & $\begin{array}{l}
  \{1-23 , 21-3 , 2-31 , 1-32\}\ or \\
  \{1-23 , 21-3 , 2-31 , 13-2\}\ or \\
  \{1-23 , 21-3 , 2-31 , 3-12\}\ or \\
  \{1-23 , 21-3 , 2-31 , 31-2\} \\
\end{array}$\\
\hline

\ref{4} & $23-1$ & $\begin{array}{l}
  \{1-23 , 2-13 , 2-31 , 1-32\}\ or \\
  \{1-23 , 2-13 , 2-31 , 13-2\}\ or \\
  \{1-23 , 2-13 , 2-31 , 3-12\}\ or \\
  \{1-23 , 2-13 , 2-31 , 31-2\} \\
\end{array}$\\
\hline

\ref{5} & $1-32$ & $\begin{array}{l}
  \{1-23 , 2-13 , 2-31 , 13-2\}\ or \\
  \{1-23 , 2-13 , 23-1 , 13-2\} \\
\end{array}$\\
\hline

\ref{5} & $1-32$ & $\begin{array}{l}
  \{1-23 , 2-13 , 13-2 , 3-12\}\ or \\
  \{1-23 , 2-13 , 13-2 , 31-2\} \\
\end{array}$\\
\hline

\ref{5} & $1-32$ & $\begin{array}{l}
  \{1-23 , 21-3 , 2-31 , 13-2\}\ or \\
  \{1-23 , 21-3 , 23-1 , 13-2\} \\
\end{array}$\\
\hline

\ref{5} & $1-32$ & $\begin{array}{l}
  \{1-23 , 21-3 , 13-2 , 3-12\}\ or \\
  \{1-23 , 21-3 , 13-2 , 31-2\} \\
\end{array}$\\
\hline

\ref{5} & $1-32$ & $\begin{array}{l}
  \{1-23 , 2-31 , 13-2 , 3-12\}\ or \\
  \{1-23 , 2-31 , 13-2 , 31-2\} \\
\end{array}$\\
\hline

\ref{5} & $1-32$ & $\begin{array}{l}
  \{1-23 , 23-1 , 13-2 , 3-12\}\ or \\
  \{1-23 , 23-1 , 13-2 , 31-2\} \\
\end{array}$\\
\hline

\ref{u1} & $12-3$ & $\begin{array}{l}
  \{1-23 , 2-13 , 2-31 , 1-32\}\ or \\
  \{1-23 , 2-13 , 2-31 , 13-2\}\ or \\
  \{1-23 , 2-13 , 23-1 , 1-32\}\ or \\
  \{1-23 , 2-13 , 23-1 , 13-2\} \\
\end{array}$\\
\hline

\ref{u1} & $12-3$ & $\begin{array}{l}
  \{1-23 , 2-13 , 2-31 , 3-12\}\ or \\
  \{1-23 , 2-13 , 2-31 , 31-2\}\ or \\
  \{1-23 , 2-13 , 23-1 , 3-12\}\ or \\
  \{1-23 , 2-13 , 23-1 , 31-2\} \\
\end{array}$\\
\hline

\ref{u1} & $12-3$ & $\begin{array}{l}
  \{1-23 , 2-13 , 1-32 , 3-12\}\ or \\
  \{1-23 , 2-13 , 1-32 , 31-2\} \\
\end{array}$\\
\hline

\ref{u2} & $12-3$ & $\begin{array}{l}
  \{1-23 , 21-3 , 23-1 , 3-12\}\ or \\
  \{1-23 , 21-3 , 23-1 , 31-2\} \\
\end{array}$\\
\hline

\ref{u2} & $12-3$ & $\begin{array}{l}
  \{1-23 , 21-3 , 2-31 , 1-32\}\ or \\
  \{1-23 , 21-3 , 23-1 , 1-32\} \\
\end{array}$\\
\hline

\ref{u2} & $12-3$ & $\begin{array}{l}
  \{1-23 , 21-3 , 2-31 , 3-12\}\ or \\
  \{1-23 , 21-3 , 2-31 , 31-2\} \\
\end{array}$\\
\hline

\end{tabular}
\end{center}
\caption{permutations avoiding five patterns}\label{5mot_2(2)}
\end{table}
\normalsize

\begin{table}
\begin{center}
\begin{tabular}{|c|c|c|}
\hline
\multicolumn{3}{|c|}{\textbf{Enumerating sequence: $\{n\}_{n\geq 1}$}}\\
\hline\hline

\emph{thanks to Proposition} &\emph{add the pattern} &\emph{to the representative} \\
\hline

\ref{3}&$3-12$&d2\\
\ref{4}&$23-1$&d9\\
\ref{3}&$3-12$&d6\\
\ref{3}&$3-12$&d25\\
\ref{2}&$21-3$&d14\\
\ref{2}&$21-3$&d17\\
\ref{5}&$1-32$&d16\\
\ref{2}&$21-3$&d20\\
\ref{4}&$23-1$&d18\\
\ref{5}&$1-32$&d21\\
\ref{u1}&$12-3$&d9\\
\ref{u1}&$12-3$&d11\\
\ref{u1}&$12-3$&d1\\
\ref{u1}&$12-3$&d2\\
\ref{u1}&$12-3$&d10\\
\ref{u1}&$12-3$&d3\\
\ref{u2}&$12-3$&d12\\
\ref{u2}&$12-3$&d4\\
\ref{u3}&$12-3$&d28\\
\ref{u3}&$12-3$&d25\\
\hline
\end{tabular}
\end{center}
\caption{permutations avoiding five patterns}\label{5mot_n}
\end{table}

\begin{table}[!h]
\begin{center}
\begin{tabular}{|c|c|c|}
\hline
\multicolumn{3}{|c|}{\textbf{Enumerating sequence: $\{F_n\}_{n\geq 1}$}}\\
\hline\hline

\emph{thanks to Proposition} &\emph{add the pattern} &\emph{to the representative} \\
\hline

\ref{u1}&$12-3$&e1\\
\ref{u1}&$12-3$&e3\\
\ref{2}&$21-3$&e4\\
\ref{4}&$1-32$&e10\\
\hline
\end{tabular}
\end{center}
\caption{permutations avoiding five patterns}\label{5mot_F_n}
\end{table}
% E' INGLOBATA IN 5MOT_SEQ_N!!!!!!!!!!!!!!!!!!!!!!!!!!!!!!!!!!!!!!!

\section{Conclusion: the cases of more than five patterns}

The approach we have followed in this work can be used to
investigate the enumeration of the permutations avoiding more than
five patterns. Really, applying the same propositions (we have
herein used) to the results about the case of the avoidance of
five patterns, one can try to solve the conjectures for the case
of six patterns. The successive cases can be examined in a similar
way.

The case of six patterns is the unique, among the remaining, which
presents some enumerating sequence not definitively constant. We
note also that all these sequences appear in the enumeration of
the case of five patterns. If $|S_n(P)|$ is required, with
$P\subseteq\mathcal M,\ |P|=6$, it should take a few minutes to
find the set $Q$ of five generalized patterns such that the
application of a certain proposition on $Q$ (among those ones
presented in this paper) leads to the set $P$ of six forbidden
patterns. So $|S_n(Q)|=|S_n(P)|$. Clearly, we are not sure that
such a set $Q$ exists since the statements in \cite{CM} are only
conjectures. Moreover, it is not sure even the fact that any
subset $P$ could be obtained by applying some proposition to some
patterns of $Q\subset P$. Nevertheless, the application of the
above mentioned propositions to the sets $Q$ of five forbidden
patterns should be confirm most of the conjectures about the case
of six patterns. This is the reason why we did not present the
analysis of this case, together with the fact that several other
tables would have appeared in these pages.

To conclude, we think that a further work about the cases of more
than six forbidden patterns does not seem to be necessary, since
many of the remaining conjectures in \cite{CM} can be easily
proved. Moreover, if $S_n(P)$ is needed, with $|P|>6$, an argument
similar to the case $|P|=6$ can be done.

%\newpage


\begin{thebibliography}{99}

\bibitem[BDPP]{BDPP} E. Barcucci, A. Del Lungo, E. Pergola, R.
Pinzani\quad\emph{ECO: A Methodology for the Enumeration of
Combinatorial Objects},\quad J. Difference Equ. Appl.\quad 5
(1999) 435-490.

\bibitem[BS]{BS} E. Babson, E. Steingr\'imsson \quad
\emph{Generalized permutation patterns and a classification of the
Mahonian statistics},\quad S\'em. Lothar. Combin.\quad 44 (2000)
B44b.

\bibitem[BFP]{BFP} A. Bernini, L. Ferrari, R. Pinzani\quad
\emph{Enumerating permutations avoiding three
Babson-Steingr\'imsson patterns},\quad Ann. Comb.\quad 9 (2005)
137-162.

\bibitem[C]{C} A. Claesson\quad \emph{Generalized pattern avoidance},
\quad Europ. J. Combin.\quad 22 (2001) 961-971.

\bibitem[CM]{CM} A. Claesson, T. Mansour\quad \emph{Enumerating
permutations avoiding a pair of Babson-Steingrímsson
patterns},\quad Ars Combin.\quad 77 (2005) 17-31.

\bibitem[EN]{EN} S. Elizalde, M. Noy\quad \emph{Consecutive patterns
in permutations},\quad Adv. in Appl. Math.\quad 30 (2003) 110-125.

\bibitem[K]{K} S. Kitaev\quad \emph{Generalized patterns in words and
permutations},\quad Ph. D. Thesis, \quad Chalmers University of
Technology and G\"oteborg Univerity (2003).

\bibitem[SS]{SS} R. Simion, F. W. Schmidt\quad\emph{Restricted
permutations}, \quad Europ. J. Combin.\quad 6 (1985) 383-406.

\end{thebibliography}
\end{document}